\documentclass[11pt,reqno]{amsart}
\usepackage{
amssymb,
amsaddr
}
\newcommand{\no}[1]{#1}

\renewcommand{\no}[1]{}  \newcommand{\upDelta}{\Delta} 

\no{\usepackage{times}\usepackage[
subscriptcorrection, slantedGreek, 
nofontinfo]{mtpro}
\renewcommand{\Delta}{\upDelta}
}
\date{\today}

\newtheorem{theorem}{Theorem}
\newtheorem{proposition}{Proposition}
\newtheorem{lemma}{Lemma}
\newtheorem{definition}{Definition}
\newtheorem{corollary}{Corollary}

\theoremstyle{remark}
\newtheorem{remark}{Remark}

\DeclareMathOperator{\supp}{supp}

\newcommand{\R}{{\bf R}}

\renewcommand{\r}[1]{(\ref{#1})}
\newcommand{\PDO}{$\Psi$DO}
\newcommand{\be}[1]{\begin{equation}\label{#1}}
\newcommand{\ee}{\end{equation}}

\renewcommand{\d}{\mathrm{d}}
\newcommand{\G}{\mathbf{G}}  \newcommand{\N}{\mathbf{N}}

\newcommand{\bo}{\partial M}
\newcommand{\Mint}{M^\text{\rm int}}

\title[Doppler transform]{The weighted Doppler transform}

\author[S. Holman]{Sean Holman}
\address[S. Holman]{Department of Mathematics \\ University of Washington \\ Seattle, WA 98195}
\author[P. Stefanov]{Plamen Stefanov}
\address[P. Stefanov]{Department of Mathematics \\ Purdue University \\ West Lafayette, IN 47907}
\thanks{Second author partly supported by NSF Grant DMS-0800428}
\date{}

\begin{document}
\maketitle
\begin{abstract}
We consider the tomography problem of recovering a covector field on a simple Riemannian manifold based on its weighted Doppler transformation over a family of curves $\Gamma$. This is a generalization of the attenuated Doppler transform. Uniqueness is proven for a generic set of weights and families of curves under a condition on the weight function. This condition is satisfied in particular if the weight function is never zero, and its derivatives along the curves in $\Gamma$ is never zero.
\end{abstract} 

\section{Introduction}

The Doppler transform of a (compactly supported) covector field $f=\{f_i\}$ in $\R^n$ is given by
\be{01}
\R^n\times S^{n-1}\ni (x,\theta)  \longmapsto \int f_j(x+t\theta)\theta^j \, \d t.
\ee
Using the Fourier transform, it can be easily seen that one can recover only $f$ up to a differential of a function $\d\phi$ so that $\phi=0$ for large $x$, see e.g., \cite{Sh-book}. In particular, such $\d\phi$ (called potential covector fields, or simply potential fields) always belong to the kernel of the transform, and this follows simply by the fundamental theorem of calculus.  The fields orthogonal to potential fields are called solenoidal and are characterized by the fact that they are divergence free. In fact, unique recovery of the solenoidal part of $f$ is possible in dimension $n = 3$ from the Doppler transform along lines parallel to only two separate planes, although stable recovery in $3$ dimensions requires the Doppler transform along lines parallel to three different planes, see \cite{Sh-Doppler}. Other partial data results and methods for recovery of the solenoidal part of $f$ from its Doppler transform can be found in \cite{Denisjuk06}, \cite{Denisjuk94}, \cite{Ramaseshan}, and \cite{Schuster3DDopp}.

The attenuated Doppler transform is defined in a similar way:
\be{02}
I_wf(x,\theta) = \int w(x+t\theta,\theta) f_j(x+t\theta)\theta^j \, \d t,
\ee
where the weight $w$ is the total attenuation along the ray $[x,x-\infty\theta]$ with an attenuation factor $\sigma(x)$
\be{03}
w(x,\theta) =  e^{-\int_{-\infty}^0\sigma(x+s\theta)\,\d s}.
\ee
It came as a surprise to find out in \cite{Buk-Kaz} (see also \cite{Natterer}) that if $\sigma$ is smooth enough and $\sigma>0$ on $\supp f$, then one can recover $f$ uniquely from $I_w$; and there are even explicit recovery formulas of the type known for the scalar attenuated transform \cite{Novikov}. 

In this work, we study a more general version of this problem. First, we work on an $n$-dimensional manifold with boundary $M$ diffeomorphic to a ball (where the support of $f$ lies), and instead of straight lines, we consider integrals over a general family of curves which we will call $\Gamma$. 
We can choose and fix a global coordinate system on $M$.
The set $\Gamma$ will have the following properties. For any $(x,\xi)\in TM\setminus \{0\}$ there exists a unique directed curve in $\Gamma$ through $x$ in the direction of $\xi$, there is a unique curve in $\Gamma$ through any given pair of points, smoothly depending on the endpoints, and $\bo$ is strictly convex w.r.t.\ $\Gamma$. In particular, the curves in $\Gamma$ solve the ODE $\ddot\gamma=G(\gamma,\dot\gamma)$ with some generator $G(x,\xi)$.  
An example are the geodesics of a  simple Riemannian metric $g$ on $M$, see the definition below. Then we study the weighted Doppler transform
\be{04}
I_wf(\gamma) = \int w(\gamma(t),\dot\gamma(t))f_j(\gamma(t))\dot\gamma^j(t)\, \d t,
\ee
and ask the following question: What condition on the weight $w(x,\xi)$ and on the family $\Gamma$ would guarantee that $I_w$ is injective? Clearly, if $w=\mbox{const.}$, then any $\d\phi$ as above is in the kernel. Moreover, if $w=\mbox{const.}$ along any curve in $\Gamma$, but maybe not a constant overall, this is still true. On the other hand, that is not true for general $w$'s and there are no obvious elements of the nullspace. One consequence of our main result is that if $w$ is never zero, and the derivative of $w$ along each $\gamma\in \Gamma$ is non-zero, then the kernel of $I_w$ is finite dimensional; generically, the kernel is trivial and there is a natural stability estimate. In particular, this applies to the attenuated Doppler transform on a simple Riemannian manifold, see section~\ref{sec4}. 
 Actually our result is more general than this. We use methods developed in \cite{SU-JAMS, FSU}. 

\section{Main results}  \label{sec2}
Let $M$ be as above; since we have assumed that $M$ is diffeomorphic to a ball we can always work in a set of global coordinates on $M$, which we will do. Let $G^i(x,\xi)$ for $i = 1, \, ... \, , n$ be a collection of smooth real-valued functions on $TM\setminus\{0\}$. We loosely think of $G(x,\xi) = (G^i(x,\xi))$ as a vector field, although it does not satisfy the transformation law required to actually make it a vector field. We will define $\Gamma$ as the curves solving $\ddot\gamma=G(\gamma,\dot\gamma)$ with given initial conditions on $\bo$. It is very convenient to extend $M$ a bit to a manifold with boundary $M_1$ with $\Mint_1\Supset M$,  so that the curves extend smoothly to $M_1$ (which can always be done if $G$ is smooth up to $\bo$). We can also assume that the so extended curves in $M_1$  still satisfy the assumptions in the introduction since those conditions are stable under small perturbations, and that $M_1$ is diffeomorphic to a slightly larger Euclidean ball. So we now have smooth functions $G^i(x,\xi)$ on $TM_1\setminus\{0\}$. We will take this as a basis of our definition and continue to describe how the curves in $\Gamma$ may be parametrized. 

Fix any smooth Riemannian metric $g$ on $M$ and extend it to $M_1$. We can always take $g$ to be Euclidean, and that is the natural choice if the curves in $\Gamma$ are not related to some metric. If there is (another) metric on $M$, and $\Gamma$ consists of geodesics of that metric, it will be more convenient to make $g$ equal to this metric, but it is not necessary. The raising and lowering of indices, and the norms below, are with respect to the metric $g$. Also, we will often refer to the unit sphere bundle $SM_1 = \{ (x,\theta) \in TM_1\; |\; |\theta|_g = 1 \}$, and define
\be{11}
\partial_{\pm}SM_1 = \left\{(x,\theta)\in \partial SM_1; \; \pm \langle\nu, \theta\rangle_g >0\right\},
\ee
where $\nu$ is the exterior unit normal to $\bo_1$, and $\partial SM_1=\{(x,\theta) \in TM_1;\; x\in \bo_1, \; |\theta|=1\}$. 

Fix a positive function $\lambda$ on $\partial_-SM_1$. Given $(x,\theta)\in \partial_-SM_1$, we define $\gamma_{x,\theta}$ as the solution of 
\be{12}
\ddot\gamma=G(\gamma,\dot\gamma), \quad \gamma(0)=x, \quad \dot\gamma(0)=\lambda(x,\theta)\theta,
\ee
defined over the largest possible  interval $[0,\tau(x,\theta)]$, $\tau(x,\theta)\ge0$. Assume that $\tau(x,\theta)$ is always finite, and also note that there is an implicit assumption in the previous sentence that the maximal interval of definition for each of the curves is closed. We now consider $\Gamma$ to be the set of curves defined by \r{12} as $(x,\theta)$ ranges over $\partial_- SM_1$, and define a smooth structure on $\Gamma$ by declaring the map $(x,\theta) \mapsto \gamma_{x,\theta}$ to be a diffeomorphism.
The generator $G$ defines a vector field $\mathbf{G}$ on $TM_1 \setminus \{0\}$  by
\be{13}
\mathbf{G} = \xi^i\frac{\partial}{\partial x^i}+G^i(x,\xi)\frac{\partial}{\partial \xi^i}.
\ee
Then $\Gamma$ consists of projections to $M_1$ of integral curves for this vector field.

We assume that $\dot\gamma$ is uniformly bounded away from zero for every $\gamma \in \Gamma$ as defined in the previous paragraph. Then the flow $\{(\gamma,\dot\gamma)\}$ of $\G$ with initial conditions as in \r{12} defines a  smooth codimension one surface $\mathcal{G}$ in $T\Mint_1$. Assume that $\mathcal{G}$ can be parametrized by $(\gamma,\dot\gamma/|\dot\gamma|_g)\in SM_1$, i.e., $\forall (y,\omega)\in SM_1$, there are  unique  $t\ge 0$, $(x,\theta)\in \partial_- SM_1$, smoothly depending on $(y,\omega)$, so that $y=\gamma(t)$, $\omega=\dot\gamma(t)/|\dot\gamma(t)|_g$ with $\gamma=\gamma_{x,\theta}$.  Then $\lambda(y,\omega)=|\dot\gamma|_g$ is a well defined smooth positive function  on $SM_1$ which can also be extended to $TM_1 \setminus \{0\}$ by homogeneity of order $-1$. That is, $\lambda(y,\theta) = \lambda(y,\theta/|\theta|_g)/|\theta|_g$. This is an appropriate extension because with it $\lambda(x,\theta) \theta$ is the velocity vector for some curve in $\Gamma$ for every $(x,\theta) \in TM$. We also extend the definition of $\gamma_{x,\theta}$ as in \r{12} but now $x$ does not need to be on $\bo$. 

Given any $(x,\theta) \in TM_1 \setminus \{0\}$ let the domain of the maximally extended curve $\gamma_{x,\theta}$ be given by $[\tau_-(x,\theta),\tau_+(x,\theta)]$. We say that $(G,\lambda)$ defines a simple system, if for any $x\in M_1$ the maps 
\[
\{(t,\theta); \; 0 < t\le\tau_+(x,\theta), \; |\theta|=1\} \ni (t,\theta) \mapsto \mathrm{exp}_{x,+}(t,\theta) =\gamma_{x,\theta}(t)\in M_1\setminus \{x\}
\]
and
\[
\{(t,\theta); \; 0 < -t\le -\tau_-(x,\theta), \; |\theta|=1\} \ni (t,\theta) \mapsto \mathrm{exp}_{x,-}(t,\theta)=\gamma_{x,\theta}(t)\in M_1\setminus \{x\}
\]
are diffeomorphisms, depending smoothly on $x$. For each $x$, we will denote the maps given above, taken together on the domain including $t=0$, by $\mathrm{exp}_x(t,\theta)$.

An example of $G$ is the generator of the geodesic flow corresponding to  $g$ on $M_1$ (or to any other metric $\tilde g$ on $M_1$). Then we can choose $\lambda=1$ (or  $\lambda=|\theta|_g/|\theta|_{\tilde g}$). The simplicity condition is fulfilled if $g$ is a simple metric \cite{Sh-book, SU-JAMS}.  Another class of examples are magnetic systems on $(M_1,g)$, where the geodesic equation is perturbed by a term representing a Lorentzian force, always normal to the velocity \cite{St-magnetic}. Then the flow preserves the energy level, so we can choose $\lambda=1$ as well. 
 
As stated in the introduction, the main object of our study will be the weighted Doppler transform
\[
I_{G,\lambda,w}:L^2(M) \to L^2(\partial_- SM_1,\mathrm{det}(g)\, \langle \theta, \nu \rangle_g\d S_x \d S_\theta)
\]
given by the formula \r{04}. Here and in what follows $L^2(M)$ refers to the space of $L^2$ covector fields. Also, the family of curves $\Gamma$ has been identified with $\partial_- SM_1$,  and the measure $\mathrm{det}(g)\, \langle \theta, \nu \rangle_g\d S_x \d S_\theta$ is chosen to provide a convenient $L^2$ structure on $\Gamma$ (other measures could also work). The measure $\mathrm{det}(g)\, \langle \theta, \nu \rangle_g \d S_x \d _\theta$ is invariantly defined but depends on the choice of $g$. The fact that the weighted Doppler transform has this mapping property can be established by considering the operator $N = I_{G,\lambda,w}^*I_{G,\lambda,w}$ which is a \PDO\ of order $-1$.
 
Before stating our main results we require one more definition.

\begin{definition}  \label{def1}
We say that $w$ satisfies the elliptic condition, if there exists an open set $U\subset \mathcal{G}$ satisfying 
\be{U}
\mbox{$\{ \xi \in T^*M \; | \; \exists \,\theta \in U$ with $\langle \xi,\theta \rangle = 0 \} = T^*M$}, \quad w|_U\not=0,
\ee
and for any $(x,\xi)\in T^*M$, there is no covector $h \in T^*_xM$ such that
\be{econd}
\G\log  w(x,\theta)= h_j \, \theta^j
\ee
for all $\theta$ normal to $\xi$ with $(x,\theta)\in U$.  
\end{definition}

\noindent Condition \r{U}, combined with the non-simplicity condition, is the necessary and sufficient condition, provided that $U$ is open, for the ellipticity of  $I_{G,\lambda,w}^* I_{G,\lambda,w}$ when acting on functions, see \cite{FSU}.

The case of dimension $n=2$ is special because then there are only two vectors $\theta$ such that $(x,\theta) \in \mathcal{G}$ is normal to any given $\xi$. If $U$ were a proper subset of $\mathcal{G}$, then for some $\xi$ the set of $(x,\theta) \in U$ with $\theta$ normal to $\xi$ would contain at most one vector, and therefore \r{econd} is trivially true for some $h$. Thus we observe that if $w$ satisfies the elliptic condition and the dimension is $2$, then it must be the case that $U = \mathcal{G}$. Taking this into consideration, in dimension $2$ we see that the elliptic condition merely states that for any $(x,\theta) \in TM \setminus \{0\}$
\[
\G \log w(x,\lambda(x,\theta)\,\theta) \neq - \frac{\lambda(x,\theta)}{\lambda(x,-\theta)} \G \log w(x,-\lambda(x,-\theta)\, \theta).
\]
In particular, if $\lambda(x,\theta)$ is even in $\theta$, which is true in the Riemannian case (where $\lambda=1$), then the condition above simply states that $\G \log w(x,\theta)$ should not be an odd function of the $\theta$ variable. 

Note that the following condition is sufficient to imply that $w$ satisfies the elliptic condition
\be{cond}
w\not=0, \quad \G w\not=0.
\ee 
This says that $w$ satisfies the elliptic condition if it does not vanish, and its derivatives along the curves in $\Gamma$ also do not vanish. In the case of the attenuated Doppler transform on a Riemannian manifold where $w$ is given by a Riemannian version of \r{03}, \r{cond} is satisfied because then $\G w=-\sigma(x)\, w(x,\theta)$, where $\sigma>0$ is the attenuation. Notice that $\sigma>0$ can depend on the direction as well, and then \r{cond} still holds. For more details, we refer to section~\ref{sec4}.

On the other hand, suppose that $h(x)$ is any smooth covector field on $M_1$. Then \r{econd} can be solved for $\log w(x,\theta)$, and therefore for $w(x,\theta)$, by integrating along curves $(\gamma,\dot \gamma) \subset TM_1$ where $\gamma \in \Gamma$. For any $\psi \in C^\infty(M_1)$ that is zero on $\partial M_1$ we then have
\[
\begin{split}
I_{G,\lambda,w}(\psi h + \d\psi)(\gamma) & = \int_\gamma w(\gamma(t),\dot \gamma(t)) (\psi(\gamma(t)) h(\gamma(t))_j + \d \psi(\gamma(t))_j) \dot\gamma(t)^j \, \d t \\
& = \int_\gamma w \psi\, h_j\, \dot\gamma(t)^j - \psi\, \G w \, \d t \\
& = \int_\gamma w \psi \, h_j\, \dot\gamma(t)^j - \psi w\, h_j\, \dot\gamma(t)^j \, \d t = 0.
\end{split}
\]
For the second equality we used integration by parts while the third follows from \r{econd}. Thus we see that some condition on $w$ is required for injectivity. In fact, with the elliptic condition we have our two main results. 

\begin{theorem}   \label{thm1}
Let w satisfy the elliptic condition, and let $G$, $\lambda$, and $w$ be real analytic on a neighborhood of the curves $(\gamma_{x,\theta},\dot \gamma_{x,\theta})$ for every $(x,\theta) \in U$ where $U$ is the set from Definition~\ref{def1}. Also suppose that $(G,\lambda)$  forms a simple system. Then  $I_{G,\lambda,w}$ is injective. Further, if $I_{G,\lambda,w}(f)(\gamma_{x,\theta}) = 0$ for all $(x,\theta) \in U$ then $f = 0$.
\end{theorem}

\begin{theorem}   \label{thm2}
Let $\lambda\in C^\infty(SM_1)$, $G\in C^\infty(TM_1 \setminus \{0\})$ form a simple system. Let $w\in C^\infty(TM_1 \setminus \{0\})$ satisfy the elliptic condition. Then the kernel of $I_{G,\lambda,w}$ is finite dimensional. Moreover, if $I_{G,\lambda,w}$ is injective then we have the following stability estimate
\be{15}
\|f^s\|_{L^2(M)} + \|\phi\|_{L^2(M)} \le C\|\mathbf{N}_{\tilde G,\tilde \lambda,\tilde w} [f,0]\|_{\mathbf{H}^1(M_1)} \quad 
\mbox{for any $f \in L^2(M)$}
\ee
and for any $(\tilde G, \tilde \lambda, \tilde w)$ in a $C^3$ neighborhood of $(G,\lambda,w)$ with some constant $C>0$. Here $f^s$ and $\phi$ are the solenoidal decomposition of $f$ on $M$ defined below, and $\mathbf{N}$ is the normal operator also defined below.
\end{theorem}

\noindent The normal operator $\mathbf{N}$ appearing on the right hand side of \r{15} is equivalent to a continuous operator applied to $I_{\tilde G,\tilde \lambda,\tilde w}$, and so this really can be called a stability estimate even though it does not provide such for $I_{\tilde G,\tilde \lambda,\tilde w}$ directly. Theorem~\ref{thm2} shows in particular that if $I_{G,\lambda,w}$ is injective for some $(G,\lambda,w) \in C^\infty$, then for every $(\tilde G,\tilde \lambda,\tilde w)$ in some $C^3$ neighborhood of $(G,\lambda,w)$, $I_{\tilde G,\tilde \lambda,\tilde w}$ is still injective. Combined with Theorem~\ref{thm1}, this gives uniqueness for an open and dense set of $(G,\lambda,w)$ with respect to the $C^3$ norm, which we call generic uniqueness. It is also true that the $C^\infty$'s in Theorem ~\ref{thm2} may be replaced by $C^k$ for $k\gg 4$.

\section{Proofs}
We start with the Hodge type decomposition of vector fields in $M$, already used in the statement of Theorem~\ref{thm2}. Given $f\in L^2(M)$, we write
\be{31}
f=f^s +\d\phi,
\ee
where $\delta f^s=0$ and $\phi\in H^1_0(M)$. Here $\delta f = \nabla^i f_i$, where $\nabla$ is the covariant derivative. 
 The function $\phi$ solves the elliptic problem
\be{32}
\Delta \phi=\delta f, \quad \phi|_{\bo}=0,
\ee
where $\Delta$ is the Laplace-Beltrami operator associated with $g$. Then $f^s\in L^2(M)$, and $\phi \in H^1_0(M)$ by elliptic regularity, see \cite{Taylor-book1}. We will also need to use the solenoidal decomposition of $f$ on $M_1$, which can be obtained by extending $f$ as zero on $M_1 \setminus M$ and replacing $\bo$ by $\bo_1$ in \r{32}. In fact, we will generally consider any function space on $M$ as the subspace of the corresponding function space on $M_1$ where the functions are are extended as zero on $M_1 \setminus M$.

Plugging \r{31} into the definition \r{04} of $I_{G,\lambda,w}$ we get
\be{33}
\begin{split}
I_{G,\lambda,w}f(\gamma) &= \int w f^s_j\dot\gamma^j\, \d t+ \int w \frac{\partial \phi}{\partial x^j}\dot\gamma^j\, \d t\\
     &= \int w f^s_j\dot\gamma^j\, \d t - \int \phi \mathbf{G}w\, \d t.
\end{split}
\ee
The problem therefore can be reformulated as follows. Given a covector field $f$ and a function $\phi$, consider the pair $[f,\phi]$. We say that $[f,\phi]$ is a solenoidal pair if $\delta f=0$.
We then study the ray transform 
\be{34}
\mathbf{I}[f,\phi](\gamma) = \int \left(w f_j\dot\gamma^j+\alpha \phi\right) \d t,
\ee
where $w(x,\xi)$ and $\alpha(x,\xi)$ are given functions on $TM$. This is actually an integral of $w(x,\xi)f_j(x)\xi^j+\alpha(x,\xi)\phi(x)$ along the maximal integral curves of $\mathbf{G}$. We are mostly interested in the special case when $\alpha = -\mathbf{G} w$, and in this case the second equality in \r{33} shows that $\mathbf{I}[f+\d \phi, \psi] = \mathbf{I}[f,\phi + \psi]$. However we can also consider the more general case when $\alpha$ may not be related to the weight $w$. 

\begin{remark} \label{ellipticalpha}
If we reformulate the elliptic condition in terms of $\alpha$, then the proofs of propositions~\ref{pr1} and~\ref{pr_par} will work in the general case when $\alpha$ does not correspond to $-\mathbf{G} w$. The appropriate reformulation is to simply replace $\mathbf{G} \log w(x,\theta)$ in \r{econd} by $\alpha(x,\theta)/w(x,\theta)$. 
\end{remark}

We consider $[f,\phi]$ as an element of the space $\mathbf{L}^2(M_1)$ with norm
\[
\|[f,\phi]\|^2_{\mathbf{L}^2(M_1)} =\int_{M_1}\bigg(   |f(x)|^2 +|\phi(x)|^2\bigg)\d v_g
\]
where $\d v_g$ is the volume form associated to $g$, and $|f|$ is the norm of the covector $f$ in the metric $g$ . Solenoidal pairs then are exactly those orthogonal to potential ones given by $[\d\psi,0]$, where $\psi\in H^1_0(M_1)$. We will denote the space of solenoidal pairs on $M_1$ (resp. $M$), which is a closed subspace of $\mathbf{L}^2(M_1)$ (resp. $\mathbf{L}^2(M)$), as $\mathbf{L}^2_s(M_1)$ (resp. $\mathbf{L}^2_s(M)$). We will also need to consider Sobolev spaces of pairs $[f,\phi]$, and will denote these, as in the statement of Theorem~\ref{thm2}, with bold letters. We view $\mathbf{I}$ as an operator $\mathbf{I} : \mathbf{L}^2(M) \to L^2(\partial_- SM_1,\mathrm{det}(g)\, \langle \theta, \nu \rangle_g\d S_x \d S_\theta)$ , and we set
\[
\mathbf{N} = \mathbf{I}^*\mathbf{I}.
\]
As mentioned above, the fact that $\mathbf{N}$ is a \PDO\ of order $-1$, which is established in Proposition~\ref{pr1}, proves that $\mathbf{I}$ does indeed map $\mathbf{L}^2(M)$ to $L^2(\partial_- SM_1,\mathrm{det}(g)\, \langle \theta, \nu \rangle_g\d S_x \d S_\theta)$. When it is necessary to distinguish between the operators $\mathbf{I}$ and $\mathbf{N}$ corresponding to different triples $(G,\lambda,w)$, or just to emphasize the dependence, we will sometimes write $\mathbf{I}_{G,\lambda,w}$ and $\mathbf{N}_{G,\lambda,w}$.

We say that $\N$ is elliptic on solenoidal pairs, if $\mbox{diag} (\Lambda \delta,0)$ and $\N$ form an elliptic system of operators, where $\Lambda$ is any parametrix of $\Delta$ in $M_1$. In other words, we want $\Lambda \delta f=0$ and $\N[f,\phi]=0$ to be an elliptic problem. 

\begin{proposition}  \label{pr1}
If $(G,\lambda)$ is a simple system, then $\mathbf{N}$ is a \PDO\ of order $-1$ in the interior of $M_1$. If $\alpha = -\mathbf{G}w$ and $w$ satisfies the elliptic condition, or $\alpha$ and $w$ satisfy the modified elliptic condition from remark~\ref{ellipticalpha}, then  $\mathbf{N}$ is elliptic on solenoidal pairs.
\end{proposition}

\begin{proof}
We follow the arguments in \cite{FSU}, where weighted integrals of functions of $x$ are considered. Note that we can think of $\xi^j$ in $w(x,\xi)f_j(x)\xi^j = \left(w(x,\xi)\xi^j\right) f_j$ as part of the weight. Then we can apply the analysis in \cite{FSU} (see the proof of Proposition~2 there) to $w\xi^{j_0}f_{j_0}$ for each $j_0$, and also to $\alpha\phi$. Summing the results we thus obtain the symbols of the corresponding \PDO s. On the other hand, $\N$ is an operator acting on pairs, unlike \cite{FSU}. The latter case is done in \cite{St-magnetic}, where $\Gamma$ consists of magnetic geodesics; see Propositions~3.9 and 4.1 there. 
Thus we get that $\N$ is an integral operator of the form
\be{35}
\N \left(\begin{matrix}f\\\phi\end{matrix}\right)  = \left( \begin{matrix}N_{11} & N_{10}\\N_{01}& N_{00}\end{matrix}\right)\left(\begin{matrix}f\\\phi\end{matrix}\right)
\ee
with  
\be{64}
\begin{split} 
(N_{11}f )^{i'}(x)&= \int_{S_xM_1} \int \lambda(x,\theta) \theta^{i'}   \bar w(x,\lambda(x,\theta) \theta) w(\gamma_{x,\theta}, \dot \gamma_{x,\theta}) 
f_i(\gamma_{x,\theta}) \dot \gamma_{x,\theta}^i J^\flat(x,\theta)  \, \d t\, \d\theta,\\
(N_{10}\phi)^{i'}(x) &= \int_{S_xM_1}  \int \lambda(x,\theta) \theta^{i'}  \bar w(x,\lambda (x,\theta) \theta) \alpha(\gamma_{x,\theta}, \dot \gamma_{x,\theta}) 
\phi(\gamma_{x,\theta})  J^\flat(x,\theta)  \, \d t\, \d\theta,\\
N_{01}f (x) &= \int_{S_xM_1}  \int  \bar \alpha(x,\lambda(x,\theta) \theta) w(\gamma_{x,\theta}, \dot \gamma_{x,\theta})  f_i(\gamma_{x,\theta}) \dot \gamma_{x,\theta}^i
 J^\flat(x,\theta)  \, \d t\, \d\theta,\\
N_{00}\phi(x) &= \int_{S_xM_1}  \int   \bar \alpha(x,\lambda(x,\theta) \theta) \alpha(\gamma_{x,\theta}, \dot \gamma_{x,\theta}) 
\phi(\gamma_{x,\theta})  J^\flat(x,\theta)  \, \d t\, \d\theta.
\end{split}
\ee
Here $J^\flat$ is an elliptic factor which is given by the Jacobian of the change of the variables $(z,\omega,t) \mapsto (x,\theta)$, where $(z,\omega)\in \partial_-SM_1$, $\gamma_{z,\omega}(t)=x$, $\dot\gamma_{z,\omega}(t)=\lambda\theta$. If $\Gamma$ consists of geodesics, or even magnetic geodesics as in \cite{St-magnetic}, then $J^\flat=1$ by the Liouville theorem. 

Following \cite{St-magnetic, FSU}, we deduce that $\N$ is a \PDO\ of order $-1$, and for the principal symbols $\sigma_p(N_{ij})$ we get
\[\begin{split} 
\sigma_p(N_{11})^{i'i}(x,\xi)&= \int_{S_xM_1}  \lambda^2(x,\theta) \theta^{i'}  \theta^i  | w(x,\lambda(x,\theta) \theta) |^2 J^\flat(x,\theta) \delta(\langle\xi,\theta\rangle) \,  \d\sigma_x(\theta),\\
\sigma_p(N_{10})^{i'}(x,\xi) &= \int_{S_xM_1}   \lambda(x,\theta) \theta^{i'}  \bar w(x,\lambda(x,\theta) \theta) \alpha(x,\lambda(x,\theta) \theta) 
 J^\flat(x,\theta)  \delta(\langle\xi,\theta\rangle) \, \d\sigma_x(\theta),\\
\sigma_p( N_{01})^i(x,\xi)  &= \int_{S_xM_1} \lambda(x,\theta) \theta^i   \bar \alpha(x,\lambda(x,\theta) \theta) w(x,\lambda(x,\theta) \theta)   
J^\flat(x,\theta) \delta(\langle\xi,\theta\rangle)  \, \d\sigma_x(\theta),\\
\sigma_p(N_{00})(x,\xi) &= \int_{S_xM_1}    |\alpha(x,\lambda(x,\theta) \theta)|^2  
 J^\flat(x,\theta)  \delta(\langle\xi,\theta\rangle)\, \d\sigma_x(\theta).
\end{split}\]
It follows that
\be{38}
\left( \sigma_p(\N)[f,\phi],[f,\phi]\right) = 
 \int_{S_xM_1} \big| \lambda(x,\theta) w(x,\lambda(x,\theta) \theta)f_j(x)\theta^j +\alpha(x,\lambda(x,\theta) \theta)\phi(x) \big|^2 J^\flat(x,\theta)  \delta(\langle\xi,\theta\rangle)\, \d\sigma_x(\theta).
\ee
Note that in the case of dimension $n = 2$, \r{38} is actually just a sum of two terms. The principal symbol of $\Lambda\delta$ is $|\xi|^{-2}\xi^i$. Fix $x$ and $\xi\not=0$. Assume now that 
\be{39}
f_i\xi^i=0
\ee 
and \r{38} vanishes. By \r{U}, the set of unit vectors $\theta$ normal to $\xi$ so that $(x,\lambda(x,\theta) \theta)\in U$, is an open subset of $\{ |\theta|=1$, $\langle \xi,\theta\rangle=0\}$. In particular, $w(x,\theta)\not=0$ for all such $\theta$.  So we have
\be{40}
f_j\lambda(\theta) \theta^j+\frac{\alpha(\lambda(\theta)\theta)}{w(\lambda(\theta)\theta)}\phi=0.
\ee
Since $x$ is fixed, we suppressed the dependence on $x$. If $\phi=0$, then \r{39} and \r{40} easily imply $f=0$ as well, and we are done. If $\phi\not=0$, then we get $\alpha(\lambda(\theta) \theta)/w(\lambda(\theta) \theta) = - \phi^{-1}f_j \lambda(\theta) \theta^j$ for $\theta$ as above. Depending on which case we consider either this last equation contradicts the modified elliptic condition given in remark~\ref{ellipticalpha}, or $\alpha = -\mathbf{G}w$ and then since $\mathbf{G}\, \mathrm{log}\, w(x,\theta) = (\mathbf{G} w) / w = -\alpha / w$ we have a contradiction of the original elliptic condition. 
\end{proof}

\begin{proposition}  \label{pr_par} 
Under the assumptions of Proposition~\ref{pr1}, the following a priori estimate holds for any solenoidal pair $[f, \phi]\in \mathbf{L}^2_s(M)$. 
\[
\|f\|_{L^2(M)} + \|\phi\|_{L^2(M)}\le  C\left(  \|\N[f,\phi]\|_{\mathbf{H}^1(M_1)} +\|[f,\phi]\|_{\mathbf{H}^{-1}(M_1)}  \right).
\]
\end{proposition}

\begin{proof}
Without loss of generality, we can assume first that $f$, $\phi$ are smooth. Also, let $M_{1/2}$ and $M_{3/4}$ be other manifolds with smooth boundary so that $M \Subset M_{1/2}^{int}\Subset M_{3/4}^{int} \Subset M_1^{int}$.

Since the pair of operators $(\N,\mathrm{diag}(\Lambda \delta,0))$ is elliptic, there exists a parametrix for the system on the interior of $M_1$. We will refer to this parametrix, which is a system of \PDO's, as $(A,B)$. Thus we have for any pair $[h,\psi] \in \mathbf{L}^2_c(M_1)$
\[
A \, \mathbf{N} [h,\psi] + B [\Lambda \delta h, 0] = [h, \psi] + K [h,\psi]
\]
where $K$ is a smoothing operator.

Now, take a cut-off function $\chi \in C^\infty_0(M_1^{int})$ that is equal to $1$ on $M_{3/4}$, and let $f^s_{M_{1}}$ be the solenoidal projection of $f$ (extended as zero to $M_1\setminus M$) in $M_{1}$. Then since $\mathrm{supp}(\delta \chi f^s_{M_1}) \subset M\setminus M_{3/4}$, by the pseudolocal property of \PDO's, $B [\Lambda \delta \chi f^s_{M_1}, 0]$ is smooth on a neighborhood of $M_{1/2}$. Thus, applying $(A,B)$ to $(\N, \mathrm{diag} (\Lambda \delta, 0)) [f^s_{M_1}, \phi]$ we have
\be{104}
\chi A\N [\chi f^s_{M_{1}},\phi] = [f^s_{M_{1}},\phi] +K[ f^s_{M_{1}},\phi] \quad \mbox{in $M_{1/2}$},
\ee
where $K$ has a smooth properly supported (in $M_1$) kernel. 

Recall that $f^s_{M_{1}}= f -\d\phi_{M_{1}}$, where 
\be{40a}
\Delta \phi_{M_{1}}=\delta f\quad \mbox{in $M_{1}$}, \qquad  \phi_{M_{1}} \in H_0^1(M_1). 
\ee
Thus the term $K[f^s_{M_{1}},\phi]$ from \r{104} can be replaced by $K[f,\phi]$ with a different $K$ with the same property.

We will write $(\Delta^D_{M_1})^{-1}$ for the solution operator of the Diricihlet realization of the Laplace-Beltrami operator $\Delta$ on $M_1$ so that with this notation $(\Delta^D_{M_1})^{-1} \delta f = \phi_{M_1}$. Now, let $\Delta^{-1}_p$ be a parametrix for the Laplacian on $M_1$ so that $((\Delta^D_{M_1})^{-1} - \Delta^{-1}_p) \delta f = K_1 f$ on $M_{1/2}$, where $K_1$ is another operator with a smooth properly supported kernel. Since $f^s_{M_{1}}= f -\d\phi_{M_{1}} = f - \d K_1f + \d \Delta_p^{-1} \delta f$ and $\chi \d = \d \chi- (\d\chi)$ we get that 
\[
A\N [\chi f^s_{M_{1}},\phi] = A\N [f ,\phi] + A \N [\d \chi K_1 f, 0] -A \N [(\d \chi) K_1f,0] + A \N [\d \chi \Delta_p^{-1} \delta f,0] - A \N [(\d \chi) \Delta_p^{-1} \delta f,0].
\] 
By \r{38}, the principal symbol of $\N$ vanishes on potential pairs, i.e., pairs of the form $[\d\psi,0]$, and so all of the terms on the right hand side except the first may be grouped together and written as a single system \PDO s of order $-1$ applied to $[f,\phi]$. Applying all of these comments to \r{104} we have
\be{40b}
A \mathbf{N}[f,\phi] = [f^s_{M_1},\phi] + K[f,\phi] \quad \mbox{in $M_{1/2}$}
\ee
where $K$ is a new \PDO\ of order $-1$.

Now, denote the projection $[f,\phi] \mapsto f$ by $\pi_1$. From \r{40b} we have
\be{41}
\d \phi_{M_{1}} = \pi_1 \big(- 
A\N [f,\phi] +  [f,\phi] +K[f,\phi]\big) \quad \mbox{in $M_{1/2}$}.
\ee
In particular,
\be{41a}
\d \phi_{M_{1}} = \pi_1 \big(- 
A\N [f,\phi] +K[f,\phi]\big) \quad \mbox{in $M_{1/2}\setminus M$}
\ee
since $[f,\phi] = 0$ outside of $M$. We integrate the equality above to get
\be{42}
\phi_{M_1} (x)= \int_0^1\langle \pi_1 \big(
A\N [f,\phi] +K[f,\phi]\big),\dot c_x \rangle \,\d s  +K_2 f     ,\quad x\in M_{1/2}\setminus M,
\ee
where $[0,1]\ni s\mapsto c_x(s)$ is any curve outside $M$ connecting $x$ and a point on $\bo_{1/2}$, smoothly depending on $x$. If $\bo_{1/2}$ is at a fixed distance to $\bo$, which we can always assume, then $c_x$ can be chosen to be the normal geodesic to $\bo$. The operator $K_2$, given by  $K_2 f(x) = \phi_{M_1}(c_x(0))$, is smoothing since $\phi_{M_1}$ is harmonic on a neighborhood of $\partial M_{1/2}$. In particular, $K_2 : L^2(M) \to H^1(M_{1/2} \setminus M)$.  

Equalities \r{41}, \r{42} yield
\be{43}
\left\|\phi_{M_1}\right\|_{H^1(M_{1/2}\setminus M)} \le C\left(  \|\N[f,\phi]\|_{\mathbf{H}^1(M_1)}) +\|[f,\phi]\|_{\mathbf{H}^{-1}(M_1)} \right).
\ee 
Using  \r{43} and the trace theorem we can also estimate $\|\phi_{M_1}\|_{H^{1/2}(\bo)}$, 
and since $f$ is solenoidal on $M$, $\phi_{M_1}$ solves $\Delta \phi_{M_{1}}= 0$ in $M$. Using the standard estimate for solutions to the Dirichlet problem we thus have
\[
\left\|\phi_{M_1}\right\|_{H^1(M)}\le C\left(  \|\N[f,\phi]\|_{\mathbf{H}^1(M_1)}) +\|[f,\phi]\|_{\mathbf{H}^{-1}(M_1)} \right).
\]
Combine this, \r{40b} and \r{41} to get the estimate in the proposition. 

\end{proof}

\begin{remark} \label{rem_Ck}
The arguments of the previous two propositions will still work for $(G,\lambda,w) \in C^k$ if $k$ is chosen large enough. When working with the \PDO\ calculus we actually only need to expand the symbols up to a certain order assuming that the remainder is still sufficiently regular.
\end{remark}

\begin{proof}[Proof of Theorem~\ref{thm2}] 

First take any solenoidal pair $[f,\phi] \in \mathrm{Ker}(\mathbf{N}_{G,\lambda,w})\cap \mathbf{L}_s^2(M)$. By Proposition~\ref{pr_par} we have that
\[
\| [f,\phi] \|_{\mathbf{L}^2(M)} \leq C \| [f,\phi] \|_{\mathbf{H}^{-1}(M_1)}.
\]
Since the inclusion $i: \mathbf{L}^2(M) \to \mathbf{H}^{-1}(M_1)$ is compact, this implies that the intersection $\mathrm{Ker}(\mathbf{N}_{G,\lambda,w}) \cap \mathbf{L}^2_s(M)$ is finite dimensional. Using \r{33} and the fact that $\mathrm{Ker}(\mathbf{I}_{G,\lambda,w}) = \mathrm{Ker}(\mathbf{N}_{G,\lambda,w})$ we see that $\mathrm{Ker}(I_{G,\lambda,w})$ is also finite dimensional, which establishes the first statement of Theorem~\ref{thm2}. We also make the observation that $\mathrm{Ker}(\mathbf{N}_{G,\lambda,w})\cap \mathbf{L}^2_s(M) \subset \mathbf{C}^{\infty}(M_1)$ ($\mathbf{C}^\infty(M_1)$ is simply the space of smooth pairs $[f,\phi]$). We can see this by using a parametrix for the pair of operators $(\mathbf{N}_{G,\lambda,w},\mathrm{diag}(\Lambda \delta,0))$

For the stability estimate we will make use of the following lemma which is similar to \cite[Prop.~V.3.1]{Taylor-book0}. See also \cite[Lemma~2]{SU-Duke}.
\begin{lemma} \label{lm_funct}
If $X$, $Y$, and $Z$ are all Banach spaces, $A:X \to Y$ is a continuous and injective linear operator, $K:X \to Z$ is a compact linear operator, and we have the estimate
\be{60}
\|x\|_X \leq C (\|Ax\|_Y + \|Kx\|_Z) \quad \forall x \in X,
\ee
then in fact we have
\[
\|x\|_X \leq \tilde C \|Ax\|_Y \quad \forall x \in X.
\]
\end{lemma}
\begin{proof}[Proof of Lemma~\ref{lm_funct}]
Suppose that the conclusion were not true and so we could find a sequence $\{x_n\}_{n=1}^\infty \subset X$ with $\|x_n\|_X = 1$ for all $n$, and such that $\|Ax_n\|_Y \to 0$ as $n \to \infty$. By compactness of $K$, there exists a subsequence $\{x_{n_k}\}_{k=1}^\infty$ so that $\{Kx_{n_k}\}_{k=1}^\infty$ converges as $k \to \infty$ and is therefore Cauchy. Thus using \r{60} and the fact $\|Ax_{n_k}\|_Y \to 0$ we see that $\{x_{n_k}\}_{k=1}^\infty$ is Cauchy. Thus $\{x_{n_k}\}$ converges to some $x \in X$. The continuity and injectivity of $A$ show that $x=0$, but this is a contradiction of the original assumption that $\|x_n\|_X = 1$ for all $n$.
\end{proof}

To apply the lemma in this situation first note that since $I_{G,\lambda,w}$ is injective and $\mathrm{Ker}(\mathbf{N}_{G,\lambda,w}) \cap \mathbf{L}_s^2(M) \subset \mathbf{C}^\infty(M_1)$, by \r{33} if $[f^s,\phi] \in \mathrm{Ker}(\mathbf{N}_{G,\lambda,w}) \cap \mathbf{L}_s^2(M)$ then $f^s + \d \phi = 0$. Since $f^s$ is solenoidal in $M$ this implies $\Delta \phi = 0$ in $M$, but also $\left . \phi \right |_{\partial M} = 0$ and so $\phi = 0$ which in turn implies $f^s = 0$. Thus $\mathbf{N}_{G,\lambda,w}$ is injective on $\mathbf{L}_s^2(M)$ and so we can take $X = \mathbf{L}_s^2(M)$, $Y = \mathbf{H}^1(M_1)$, $Z = \mathbf{H}^{-1}(M_1)$, $A = \mathbf{N}_{G,\lambda,w}$, and $K$ the inclusion map from $X$ to $Z$. Using the estimate from Proposition~\ref{pr_par}, we establish
\be{61}
\|[f^s,\phi]\|_{\mathbf{L}^2(M)} \leq C \|\mathbf{N}_{G,\lambda,w} [f^s,\phi]\|_{\mathbf{H}^1(M_1)} \quad \forall [f^s,\phi]\in \mathbf{L}^2_s(M)
\ee
for some constant $C>0$. If $f \in L^2(M)$, we may now take its solenoidal decomposition $f = f^s + \d \phi$ on $M$. Applying \r{61} to $[f^s,\phi]$ yields
\be{62}
\|f^s\|_{L^2(M)} + \|\phi\|_{L^2(M)} \leq  C \|\mathbf{N}_{G,\lambda,w} [f^s,\phi]\|_{\mathbf{H}^1(M_1)}.
\ee
If we note that by \r{33}, $\mathbf{N}_{G,\lambda,w}[f,0] = \mathbf{N}_{G,\lambda,w}[f^s,\phi]$ we can see that this is the stability estimate \r{15} in the case $(\tilde G,\tilde \lambda,\tilde w) = (G,\lambda,w)$.

To complete the proof we argue that for small enough perturbations $(\tilde G,\tilde \lambda,\tilde w)$ of $(G,\lambda,w)$ in $C^3$ the estimate \r{62} still holds with a constant $C$ independent of $(\tilde G, \tilde \lambda, \tilde w)$. Note that in the case that $(G,\lambda,w)$ are only $C^3$, we may still define $I_{G,\lambda,w}$ and the normal operator $\mathbf{N}_{G,\lambda,w}$ which is still given by \r{35}. We will follow \cite{FSU} and establish the following proposition.
\begin{proposition} \label{pr_C2}
If $(G,\lambda,w) \in C^\infty$ is as in the statement of Theorem~\ref{thm2}, and $(\tilde G,\tilde \lambda,\tilde w) \in C^3$ with $\|(\tilde G,\tilde \lambda,\tilde w) - (G,\lambda,w)\|_{C^3} \leq \delta$, then for all $\delta$ sufficiently small $(\tilde G,\tilde \lambda)$ is still a simple system (where the diffeomorphisms are only $C^3$), and there exists a constant $C >0$, which depends only on $(G,\lambda,w)$ and the dimension $n$ such that
\be{63}
\| (\mathbf{N}_{G, \lambda, w} - \mathbf{N}_{\tilde G,\tilde \lambda,\tilde w}) [f,\phi] \|_{\mathbf{H}^1(M_1)} \leq C \delta (\|f\|_{L^2(M)} + \|\phi\|_{L^2(M)}) \quad \forall f,\phi \in L^2(M).
\ee
\end{proposition}
Before proving Proposition~\ref{pr_C2} we note that the proof of Theorem~\ref{thm2} follows if we take a new constant $C$ which is the product of the two constants from \r{62} and \r{63} and then make $\delta < 1/2C$. Using both \r{62} and \r{63} then gives
\[\begin{split}
\|f^s\|_{L^2(M)} + \|\phi\|_{L^2(M)} \leq & \; C \left ( \|\mathbf{N}_{\tilde G,\tilde \lambda,\tilde w} [f^s,\phi] \|_{\mathbf{H}^1(M_1)} + \| (\mathbf{N}_{G, \lambda, w} - \mathbf{N}_{\tilde G,\tilde \lambda,\tilde w}) [f^s,\phi] \|_{\mathbf{H}^1(M_1)} \right ) \\ 
\leq &\; C\;  \|\mathbf{N}_{\tilde G,\tilde \lambda,\tilde w} [f^s,\phi] \|_{\mathbf{H}^1(M_1)} + \frac{1}{2} (\|f^s\|_{L^2(M)} + \|\phi\|_{L^2(M)} )
\end{split}\]
and so
\[
\|f^s\|_{L^2(M)} + \|\phi\|_{L^2(M)} \leq 2C \; \|\mathbf{N}_{\tilde G,\tilde \lambda,\tilde w} [f,0] \|_{\mathbf{H}^1(M_1)} \quad \forall f \in L^2(M)
\]
and for all $(\tilde G, \tilde \lambda, \tilde w)$ in the $\delta$ neighborhood of $(G,\lambda,w)$ with respect to the $C^3$ norm. This completes the proof of Theorem~\ref{thm2}.
\end{proof}
\begin{proof}[Proof of Proposition~\ref{pr_C2}]
The proof will consist of a careful comparison of the integral kernels corresponding to $\mathbf{N}_{G,\lambda,w}$ and $\mathbf{N}_{\tilde G,\tilde \lambda, \tilde w}$ given by \r{64}. We begin by considering the two maps
\be{69}
F_x,\tilde F_x:(t,\theta) \mapsto (r,\omega) = \left (\mathrm{sign}(t)\, |\mathrm{exp}_x(t,\theta) - x|, \mathrm{sign}(t)\, \frac{\mathrm{exp}_x(t,\theta) - x}{|\mathrm{exp}_x(t,\theta) - x|} \right ).
\ee
Here $\mathrm{exp}_x$ may correspond to either $(G,\lambda)$ or $(\tilde G, \tilde \lambda)$, and the corresponding maps given by \r{69} are denoted by $F_x$ and $\tilde F_x$ respectively. We extend both $(G,\lambda)$ and $(\tilde G, \tilde \lambda)$ sufficiently far, possibly beyond $M_1$, so that these two maps have the same domain. The $(t,\theta)$ coordinates are loosely speaking polar coordinates with respect to the family of curves generated by $G$ (or $\tilde G$) while $(r,\omega)$ are actual (Euclidean) polar coordinates centered at $x$. 

Our first order of business will be to estimate $F_x^{-1} - \tilde F_x^{-1}$ in the $C^2_{x,r,\omega}$ norm on the set of $(x,r,\omega)$ such that $x + r\omega \in M_1$. We will accomplish this in three steps. First we will estimate $(F_x - \tilde F_x)$, and then using this we estimate $(F_x^{-1} \circ \tilde F_x - Id)$. Finally, we precompose by $\tilde F^{-1}_x$ to obtain the desired estimate.

For the first step we use the following standard estimate from ODE theory which is proven in \cite{CL} and also applied in the same manner in \cite{FSU}.
\begin{lemma} \label{lm_ODE}
Let $x$ and $\tilde x$ solve the ODE systems
\[
x' = F(t,x), \quad \tilde x' = \tilde F(t,\tilde x),
\]
where $F$, $\tilde F$ are continuous functions from $[0,T] \times U$ to a Banach space $\mathcal{B}$, where $U \subset \mathcal{B}$ is open. Let $F$ be Lipschitz w.r.t. $x$ with a Lipschitz constant $k>0$. Assume that
\[
\|F(t,x) - \tilde F(t,x)\| \leq \delta, \quad \forall t \in [0,T], \; \forall x \in U,
\]
and that $x(t)$, $\tilde x(t)$ stay in $U$ for $0\leq t \leq T$. Then for $0 \leq t \leq T$
\[
\|x(t) - \tilde x(t) \| \leq e^{kt} \|x(0) - \tilde x(0)\| + \frac{\delta}{k}\left ( e^{kt} - 1 \right ).
\]
\end{lemma}

\noindent We apply the Lemma~\ref{lm_ODE} to the systems
\be{100}
\begin{array}{ll}
\gamma'(t) = \xi(t) & \tilde{\gamma}'(t) = \tilde{\xi}(t) \\
\xi'(t) = G(\gamma(t),\xi(t)) & \tilde{\xi}'(t) = \tilde{G}(\tilde{\gamma}(t),\tilde{\xi}(t))
\end{array}
\ee
considered with the same initial conditions $\gamma(0) = \tilde{\gamma}(0) = x$ and $\xi(0) = \tilde{\xi}(0) = \lambda(x,\theta)\, \theta$ to obtain
\be{72}
\| \gamma_{x,\theta} - \tilde \gamma_{x,\theta}\|_{C^3_{x,t,\theta}} + \| \dot{\gamma}_{x,\theta} - \dot{\tilde{ \gamma}}_{x,\theta}\|_{C^3_{x,t,\theta}} = O( \delta).
\ee
To get the estimates for the $x$ and $\theta$ derivatives we first differentiate the systems \r{100} with respect to the initial conditions, and then apply Lemma~\ref{lm_ODE} as well as the hypotheses about $(G,\lambda)$ and $(\tilde G,\tilde \lambda)$. For the $t$ derivatives we appeal directly to the original system \r{100} and use the assumptions on $G$ and $\tilde G$. For the mixed derivatives we use a combination of these techniques.

At this point note that \r{72} proves the statement that $(\tilde G, \tilde \lambda)$ is still a simple system, and together with the fact that we may rewrite \r{69} in terms of $\dot \gamma_{x,\theta}$ (resp. $\dot{ \tilde{ \gamma}}_{x,\theta}$) it also shows that
\be{70}
\|F_x(t,\theta) - \tilde F_x(t,\theta)\|_{C^3_{x,t,\theta}} = O(\delta).
\ee
Now note that $(F_x^{-1} \circ \tilde F_x - Id) = (F_x^{-1} \circ \tilde F_x - F_x^{-1} \circ F_x)$, and so, working in some appropriate set of local coordinates for $\omega \in \mathbb{S}^{n-1}$, we have
\[
(F_x^{-1} \circ \tilde F_x - Id)(t,\theta) = \left ( \int_{0}^1 DF_x^{-1}(s\tilde F_x(t,\theta) + (1-s) F_x(t,\theta)) \, \d s \right ) \cdot (\tilde F_x(t,\theta) - F_x(t,\theta)).
\]
Taking derivatives of his last equation we see that $\|F_x^{-1} \circ \tilde F_x - Id\|_{C^3_{x,t,\theta}}$ can be bounded in terms of \r{70}, and $\|DF_x^{-1}\|_{C^3_{x,r,\omega}}$. By \r{70}, $\|\tilde F_x\|_{C^3_{x,t,\theta}}$ and $\|D\tilde{F}_x^{-1}\|_{C^2_{x,r,\omega}}$ are uniformly bounded for $(\tilde G, \tilde \lambda)$ in a sufficiently small $C^3$ neighborhood of $(G,\lambda)$. Thus, if we precompose $F_x^{-1} \circ \tilde F_x - Id$ with $\tilde F_x^{-1}$ we see that
\be{71}
\|F_x^{-1}(r,\omega) - \tilde F_x^{-1}(r,\omega)\|_{C^3_{x,r,\omega}} = O(\delta).
\ee

Now we use the maps $F_x^{-1}$ and $\tilde F_x^{-1}$ to change variables in \r{64} for $\mathbf{N} = \mathbf{N}_{G,\lambda,w}$ and $\tilde{\mathbf{N}} = \mathbf{N}_{\tilde G,\tilde \lambda,\tilde w}$ respectively. We will call the Jacobian determinants of these changes, which depend on $x$, $J(r,\omega,x)$ and $\tilde{J}(r,\omega,x)$. By \r{71} we have that $\|J(x,r,\omega) - \tilde{J}(x,r,\omega)\|_{C^2_{x,r,\omega}} = O(\delta)$. We also introduce a cutoff function $\chi \in C_c^\infty(M_1)$ equal to $1$ on $M$, and since the integrands in \r{64} are zero outside of $M$ we may multiply them by $\chi(x+r\omega)$ without changing the values of the integrals. Thus the first equation for $\bf{N}$ from \r{64} becomes
\be{65}
(N_{11}f )^{i'}(x) = \int_{S_xM_1} \int_{\mathbb{R}} \chi(x+r\omega) A_{11}^{ii'} (x,r,\omega) f_i(x+r\omega) \, \d r \d\omega.
\ee
where
\[
A_{11}^{ii'} (x,r,\omega) = \left ( \lambda (x,  \theta) \bar w (x,\lambda(x,\theta) \theta) w (\gamma_{x,\theta}(t), \dot \gamma_{x,\theta}(t))  \theta^i \dot \gamma_{x,\theta}(t)^{i'}
J^\flat (x,\theta) J(x,r,\omega) \right |_{(t,\theta) = F_x^{-1}(r,\omega)}.
\]
We have the same equation for $\tilde{\bf{N}}$ when tildes are added in appropriate places. Also, we can find similar equations for $N_{10}$, $N_{01}$, and $N_{00}$, and the following analysis proceeds in the same manner for those operators. We will only examine $N_{11}$ in detail.  

By changing variables $(r,\omega) \mapsto (-r,-\omega)$ we see that if $A_{11}^{ii'}$ is odd with respect to $(r,\omega)$, then \r{65} is always zero. Thus we may replace $\chi A_{11}^{ii'}$ by $(A_e)_{11}^{ii'} = \chi(x+r\omega) (A_{11}^{ii'}(x,r,\omega) + A_{11}^{ii'}(x,-r,-\omega))/2$ which is even with respect to $(r,\omega)$. Having done this we take the linear approximation in $r$ for $(A_e)_{11}^{ii'}$
\be{66}
(A_e)_{11}^{ii'}(x,r,\omega) = (A_e)_{11}^{ii'}(x,0,\omega) + r R^{ii'}(x,r,\omega).
\ee
Provided that $A_{11}^{ii'}\in C^2_{x,r,\omega}$ we see that $R^{ii'} \in C^1_{x,r,\omega}$. Changing to Cartesian coordinates $y = x + r\omega$ in \r{65} yields
\[
\begin{split}
(N_{11}f)^{i'} (x) =\;
& 2\int_{\mathbb{R}^n} (A_e)_{11}^{ii'}(x,0,\frac{y-x}{|y-x|}) \frac{f_i(y)}{|y-x|^{n-1}} \d y \\
& + 2\int_{\mathbb{R}^n} R^{ii'}(x,|y-x|,\frac{y-x}{|y-x|}) \frac{f_i(y)}{|y-x|^{n-2}} \d y.
\end{split}
\]
In the same manner as above we find kernels $\tilde A_e$ and $\tilde R$ associated to $(\tilde G,\tilde \lambda,\tilde w)$. Thus, comparing the two operators we have
\be{67}
\begin{split}
((N_{11}- \tilde N_{11})f)^{i'} (x) = \;
& 2\int_{\mathbb{R}^n} \left ((A_e)_{11}^{ii'}(x,0,\frac{y-x}{|y-x|}) - (\tilde A_e)_{11}^{ii'}(x,0,\frac{y-x}{|y-x|}) \right ) \frac{f_i(y)}{|y-x|^{n-1}} \d y \\
& + 2\int_{\mathbb{R}^n} \left (R^{ii'}(x,|y-x|,\frac{y-x}{|y-x|}) - \tilde R^{ii'}(x,|y-x|,\frac{y-x}{|y-x|}) \right ) \frac{f_i(y)}{|y-x|^{n-2}} \d y.
\end{split}
\ee
Let $A(x,\omega)$ and $R(x,r,\omega)$ denote respectively the portions of the first and second integral kernels from \r{67} contained in parentheses and assume for the moment that
\be{68}
\|A(x,\omega)\|_{C^1_{x,\omega}} + \|R(x,r,\omega)\|_{C^1_{x,r,\omega}} < C \delta
\ee
for a constant C. Assuming $f \in L^2(M)$, we may estimate easily the second integral in \r{67} and its first derivatives in $L^2(M_1)$ by $C \delta$ since the kernels of both will have only integrable singularities. Similarly, the first integral in \r{67} may be estimated in $L^2(M_1)$ by $C \delta$. However, the derivatives of the first integral pose a problem since when the kernel is differentiated its singularity is no longer integrable. However we may still use the Calderon-Zygmund theorem to estimate the first integral with a differentiated kernel in $L^2$, and by \cite[Theorem XI.11.1]{MikhlinP} this is the derivative of the integral. Therefore all that remains to prove the proposition is to prove the estimate \r{68}.

To prove \r{68} we first note, from the definition of $A(x,\omega)$ and $R(x,r,\omega)$, that it is sufficient to prove
\[
\|A^{ii'}_{11}(x,r,\omega) - \tilde A^{ii'}_{11}(x,r,\omega)\|_{C^2_{x,r,\omega}} = O(\delta).
\]
We require $C^2$ since $R(x,r,\omega)$ already involves one derivative. It is thus sufficient to show that the differences of each of the corresponding terms in the definitions of $A^{ii'}_{11}(x,r,\omega)$ and $\tilde A^{ii'}_{11}(x,r,\omega)$ (see below \r{65}) are $O(\delta)$ in $C^2_{x,r,\omega}$. Indeed, this has already been shown for $(J(x,r,\omega) - \tilde{J}(x,r,\omega))$. The other terms all involve $(t(x,r,\omega),\theta(x,r,\omega))$, and so we can estimate them in $C^2_{x,r,\omega}$ by using \r{71} together with an estimate of the difference in $C^2_{x,t,\theta}$, and a $C^{2,\infty}_{x,t,\theta}$ bound for $(G,\lambda,w)$ (the bound must be in $C^{3,\infty}_{x,t,\theta}$ when we are looking at $\alpha = -\G w$). Here and in what follows the $C^{k,\infty}_{x,t,\theta}$ norm is the usual $C^k_{x,t,\theta}$ norm plus the Lipschitz norm of all $k$th derivatives.

As an example of this method we estimate $\lambda(x,\theta(x,r,\omega)) - \tilde \lambda(x,\tilde \theta(x,r,\omega))$ in $C^2_{x,r,\omega}$ where $\theta$ and $\tilde \theta$ represent the projections onto the $\theta$ component of $F_x^{-1}$ and $\tilde F_x^{-1}$ respectively. First we have
\[
\begin{split}
\| \lambda(x,\theta(x,r,\omega)) - \tilde \lambda(x,\tilde \theta(x,r,\omega)) \|_{C^2_{x,r,\omega}} \leq \; &
\| \lambda(x,\theta(x,r,\omega)) - \lambda(x,\tilde \theta(x,r,\omega)) \|_{C^2_{x,r,\omega}} \\ &
+ \| \lambda(x,\tilde \theta(x,r,\omega)) - \tilde \lambda(x,\tilde \theta(x,r,\omega)) \|_{C^2_{x,r,\omega}}.
\end{split}
\]
We now see that the first term above can be bounded by \r{71} and a $C^{2,\infty}_{x,t,\theta}$ bound on $\lambda$. The second term is $O(\delta)$ by the hypotheses, and so in fact the entire right hand side is $O(\delta)$. The same method, also using \r{72}, works for the difference between the $w$ and $\tilde w$ terms. For the difference between $\alpha$ and $\tilde \alpha$, which is not required for $N_{11}$ but does arise when we consider $N_{01}$ and $N_{00}$, the same method works, but because $\alpha$ already involves one derivative of $w$ we need that $\|w - \tilde w\|_{C^3_{x,t,\theta}} = O(\delta)$ and require a bound on $w$ in $C^{3,\infty}_{x,t,\theta}$. Finally, for $(J^\flat - \tilde J^\flat)$ we use \r{72}, \r{71}, together with the smoothness of the boundary of $M_1$. This completes the proof of Proposition~\ref{pr_C2}.

\end{proof}

We next move on to the proof of Theorem~\ref{thm1}. Here we will not make use of the normal operator $\N$, but rather analyze $I_{G,\lambda,w}$ itself in a neighborhood of a single curve of $\Gamma$. Our method is analytic microlocal in nature and relies on the calculus of analytic wavefront sets which we will denote by $\mathrm{WF_A}([f,\phi])$ (see \cite{Sj-A}). The same technique has previously been applied to other problems: see \cite{FSU}, \cite{SU-AJM}. The main step is to establish the following proposition.

\begin{proposition} \label{pr3}
Let $G$, $\lambda$, $w$, and $U$ be as in the hypotheses of Theorem ~\ref{thm1}. Suppose $[f^s,\phi] \in \mathbf{L}^2(M_1)$ is a pair such that $f^s$ is solenoidal in the interior of $M$, $\mathbf{I}[f^s,\phi](\gamma_{x,\theta}) = 0$ for every $(x,\theta) \in U$, and $\mathrm{supp}([f^s,\phi]) \subset M$. Then $[f^s,\phi]$ is analytic in the interior of M.
\end{proposition}

\begin{proof}
This proof will follow closely that of the corresponding results in \cite{FSU} and \cite{SU-AJM}. Take any $(x_0,\xi^0) \in T^* M$ with $x_0$ in the interior of $M$. We aim to show that
\[
(x_0, \xi^0) \notin \mathrm{WF_A}([f^s,\phi]).
\]
The result then follows since the projection of $\mathrm{WF_A}([f^s,\phi])$ to $M$ is exactly the set of points where $[f^s,\phi]$ is not analytic.

By the elliptic condition, there exists $(x_0,\theta_0) \in U$ normal to $(x_0,\xi^0)$. Since $(x_0,\theta_0)$ is  in $\mathcal{G}$,  there must be a point $(p_0,\omega_0) \in \partial_-SM_1$ such that $\gamma_{p_0,\omega_0}$ passes through $(x_0,\theta_0)$. We will now define coordinates in a neighborhood of $\gamma_{p_0,\omega_0}$. First, let us take some coordinates $(\omega_1,\, ...\, , \omega_{n-1})$ on $S_{p_0} M_1$ centered at $\omega_0$. Then by the simplicity assumption $(\omega_1, \,...\, , \omega_{n-1},t) \mapsto \mathrm{exp}_{p_0}(t,\omega)$ is the inverse of a coordinate map on a neighborhood of $\gamma_{p_0,\omega_0}$. We will write $x' = (x^1, \, ... \,, x^{n-1}) = (\omega^1, \, ...\, , \omega^{n-1})$ and $x^n = t$. Next we translate the coordinates so that $x_0 = 0$ and now assume that they are defined on the set
\[
V = \{ |x'| < \epsilon,\; l_- < x^n < l_+ \} \subset M_1
\]
where $l_-$ and $l_+$ are such that $(x',l_-)$ and $(x',l_+)$ lie outside of $M$ for $|x'| <\epsilon$. On a small enough neighborhood we may assume that these coordinates are analytic because $G$ and $\lambda$ are analytic on a neighborhood of $\gamma_{p_0,\omega_0}$.  Also note that in these coordinates $(0,(0,1)) = (x_0,\theta_0)$.

Additionally, if we fix any $\theta' = (\theta^1, \, ... \, , \theta^{n-1})$, then for such $\theta'$ sufficiently small the map $(z', t) \mapsto \mathrm{exp}_{(z',0)}(t,(\theta',1))$ defines another analytic coordinate system, which depends analytically on $\theta'$, on some subset of $V$. In particular, if we restrict to $|z'| < 7 \epsilon/8$, then for $\theta'$ small enough the curves $\gamma_{(z',0),(\theta',1)}$ reach points outside of $M$ while still inside $V$. Thus since $U \subset \mathcal{G}$ is open, by taking $|\theta'| \ll 1$ we can ensure that $\mathbf{I}[f^s,\phi](\gamma_{(z',0),(\theta',1)}) = 0$.

Now we introduce a sequence of cutoff functions $\chi_N \in C_c^\infty(\mathbb{R}^{n-1})$ such that $\mathrm{supp}(\chi_N) \subset \{|z'| < 3 \epsilon/4\}$, $\chi_N (z') = 1$ for $|z'|<\epsilon/2$, and 
\be{45}
|\partial^\tau \chi_N (z')| < (CN)^{|\tau|} \quad , \quad \forall\, z', |\tau| < N.
\ee
It is possible to construct such cutoff functions, see \cite[Chapter V, Lemma 1.1]{Treves}.

For all $\xi$ in a complex neighborhood of $\xi^0$ and $h > 0$ we multiply $\mathbf{I}[f^s,\phi](\gamma_{(z',0),(\theta',1)}) = 0$ by $\chi_N(z') e^{\frac{i}{h} (z' \cdot \xi)}$ and integrate in $z'$ to obtain
\be{46}
\iint e^{\frac{i}{h} (z'\cdot\xi)} \chi_N(z') (w f_j^s \dot \gamma_{(z',0),(\theta',1)}^j + \alpha \phi)\, \d t \, \d z' = 0.
\ee
As observed above, $(z',t)$ provide coordinates on a subset of $V$, and outside this subset the integrand in the following formula is $0$. Thus we may change coordinates in \r{46} to obtain
\be{47}
\int e^{\frac{i}{h} (z'(x,\theta') \cdot \xi)} (a_N(x,\theta')\, f^s_j \, b^j(x,\theta') + c_N(x,\theta')\, \phi) \, \d x = 0
\ee
where for sufficiently small $\theta'$ we have that $a_N$ and $c_N$ are analytic and independent of $N$ on a neighborhood of $\gamma_{p_0,\omega_0}$, and satisfy \r{45}, with a new constant $C$, in all variables everywhere. Also, $a_N(0,\theta') = w(0,\lambda(\theta',1) (\theta',1))$ and $c_N(0,\theta') = \alpha(0,\lambda(\theta',1)(\theta',1))$. Finally, the components $b^j(x,\theta')$ of the vector field $b$ are analytic everywhere and we have $b(0,\theta') = \lambda(0,(\theta',1)) (\theta',1)$.

Now we will choose $\theta'$ based on $\xi$. For convenience we will write $\theta = (\theta',1)$. First, without loss of generality we may assume that $\xi^0 = \d x^{n-1}$. Now, following \cite{FSU} and \cite{SU-AJM} we are able to define an analytic mapping $\theta'(\xi)$, and corresponding $\theta(\xi)$, for $\xi$ in a neighborhood of $\xi^0 = \d x^{n-1}$ with the following properties.
\[
\theta(\xi)\cdot \xi = 0, \quad \theta^n(\xi) = 1, \quad \theta(\xi^0) = \theta_0,
\]
and
\[
\frac{\partial \theta}{\partial \xi_\nu} = \frac{\partial}{\partial x^\nu}, \quad \nu = 1, \,...\, , n-2, \quad \frac{\partial \theta}{\partial \xi_{n-1}}(\xi_0) = 0, \quad \frac{\partial \theta}{\partial \xi_n}(\xi_0) = -\frac{\partial}{\partial x^{n-1}}.
\]
Note that in order to define a function with all of these properties we need that $\theta_0 \cdot \xi^0 = 0$. Replacing $\theta'$ by $\theta'(\xi)$ in \r{47} we arrive at
\be{48}
\int e^{\frac{i}{h} \psi(x,\xi)} (\tilde a_N(x,\xi)\, f^s_j \, \tilde b^j(x,\xi) + \tilde c_N(x,\xi)\, \phi)\, \d x = 0
\ee
where $\tilde a_N$, $\tilde b^j$, and $\tilde c_N$ have the same properties as $a_N$, $b^j$, and $c_N$ respectively, and
\[
\psi(x,\xi) = z'(x,\theta'(\xi)) \cdot \xi.
\]
Our next step will be to apply the method of complex stationary phase (see \cite{Sj-A}), and so we must analyze the critical points of this phase function. In fact, $\psi$ is precisely the same as the phase considered in \cite{FSU}, and so we take two results from that paper. Firstly, we have
\be{52}
\psi_{x\xi}(0,\xi) = \mathrm{Id},
\ee
and second we have the following lemma.
\begin{lemma}\label{lm1}
There exists $\delta >0$ such that
\[
\psi_\xi(x,\xi) \neq \psi_\xi(y,\xi) \quad for \; x \neq y,
\]
for $x \in V$, $|y| < \delta$, $|\xi - \xi^0|<\delta$, $\xi$ complex.
\end{lemma}
It is in the proof of this lemma that $\theta(\xi) \cdot \xi = 0$ is used.

Now for $y$ and $\delta$ as in lemma~\ref{lm1} and $\eta \in B(\xi^0,\delta/2)$, we multiply \r{48} by
\[
\tilde \chi(\eta - \xi) e^{\frac{i}{h}(\frac{i}{2}(\xi - \eta)^2 - \psi(y,\xi))},
\]
where $\tilde \chi$ is the characteristic function of $B(0,\delta/2)$, and then integrate in $\xi$. This yields
\be{49}
\iint e^{\frac{i}{h} \Psi(y,x,\eta,\xi)} (A_N(x,\xi,\eta)\, f^s_j \, \tilde b^j(x,\xi) + C_N(x,\xi,\eta) \, \phi)\, \d x\, \d \xi = 0
\ee
where
\[
\Psi(x,y,\eta,\xi) = -\psi(y,\xi) + \psi(x,\xi) + \frac{i}{2} (\xi - \eta)^2.
\]
$A_N$ and $C_N$ simply incorporate the dependence on $\eta$ into $\tilde a_N$ and $\tilde c_N$. Let us first consider the portion of the integral in \r{49} where $|x-y| > \delta/C^1$ for some constant $C^1 >1$. By lemma~\ref{lm1} there is no real critical point for the function $\xi \mapsto \Psi(x,y,\eta,\xi)$ on this set, and thus we can estimate that portion of the integral as follows.
\be{50}
\begin{split}
\left | \iint_{|x-y|>\delta/C^1} e^{\frac{i}{h} \Psi(y,x,\eta,\xi)} (A_N(x,\xi,\eta)\, f^s_j \, \tilde b^j(x,\xi) + C_N(x,\xi,\eta) \, \phi)\, \d x\, \d \xi \right | \\ < C^2 \left (C^2 N h \right)^N + C^2 N e^{-1/(h C^2)}
\end{split}
\ee
for a new constant $C^2$. To obtain \r{50} we integrate by parts $N$ times in $\xi$ using the facts that $|\Psi_\xi|$ is bounded below on the domain of integration,
\[
e^{\frac{i}{h} \Psi(x,y,\eta,\xi)} = h \frac{\overline \Psi_\xi \cdot \partial_\xi}{i |\Psi_\xi|^2} e^{\frac{i}{h} \Psi(x,y,\eta,\xi)},
\]
and on the boundary of the domain of integration in $\xi$ (where $|\xi - \eta| = \delta/2$), $e^{\frac{i}{h} \Psi(x,y,\eta,\xi)}$ is exponentially small in $1/h$. We also made use of the fact that the amplitude satisfies an estimate like \r{45} in $\xi$.

Now let us take a closer look at the critical points of $\xi \mapsto \Psi(x,y,\eta,\xi)$. We easily see that for $x = y$ there is a unique critical for this function, namely $\xi_c = \eta$. Combining this with \r{52} and the implicit function theorem we see that for $|x-y|$ sufficiently small there is still a unique complex critical point $\xi_c$, depending analytically on $x$, $y$, and $\eta$. Further, $\mathcal{I}\Psi_{\xi\xi}(x,x,\eta,\xi) = \mathrm{Id} > 0$ and so it is still true that $\mathcal{I}\Psi_{\xi\xi}(x,y,\eta,\xi_c)>0$ if $|x-y|<\delta/C^1$ with $C^1$ sufficiently large (this is required to apply the complex method of stationary phase below). Making $C^1$ possibly larger still we can ensure that $A_N$ and $C_N$ are analytic and independent of $N$ on $|x-y| < \delta/C^1$. Thus we may apply the complex method of stationary phase (see \cite{Sj-A}, Theorem 2.8 and the remark after the theorem) to the rest of the integral in \r{49} to obtain
\be{51}
\begin{split}
\iint_{|x-y|\leq\delta/C^1} e^{\frac{i}{h} \Psi(x,y,\eta,\xi)} (A_N(x,\xi,\eta)\, f^s_j \, \tilde b^j(x,\xi) + C_N(x,\xi,\eta) \, \phi)\, \d x\, \d \xi \\ 
= \int_{|x-y|\leq\delta/C^1} e^{\frac{i}{h} \Psi(x,y,\eta,\xi_c)} (A(x,\beta;h) \, f^s_j \, b^j(x,\beta) + C(x,\beta;h) \, \phi) \, \d x+ \mathcal{O}(e^{-\frac{1}{hC^3}})
\end{split}
\ee
where we are using the notation $\beta = (y,\eta)$, and $b^j(x,\beta)$ is a new vector valued function. We will also write
\[
\Phi(x,\beta) = \Psi(x,y,\eta,\xi_c).
\]
Expanding $\Phi$ about $y=x$ we have
\[
\Phi(x,\beta) = -\psi_x(x,\eta) (y-x) + (y-x)^t \left ( -\psi_{xx}(x,\eta) + i \psi_{x\xi}^2(x,\eta) \right ) (y-x) + \mathcal{O}(|y-x|^3),
\]
and so if $\delta$ is sufficiently small, by \r{52} we have that $\mathcal{I}\Phi(x,\beta) > C^4 |y-x|^2$ for a positive constant $C^4$ on the domain of integration in \r{51} when $\beta$ is real. Now, combining \r{49}, \r{50}, and \r{51} we have
\be{53}
\begin{split}
\int_{|x-y|\leq\delta/C^1} e^{\frac{i}{h} \Phi(x,\beta)} (A(x,\beta;h) \, f^s_j \, b^j(x,\beta) + C(x,\beta;h) \, \phi) \, \d x \\ 
= \mathcal{O}\left (N e^{-\frac{1}{hC^3}} + \left (C^2 N h \right )^N \right ),
\end{split}
\ee
and by choosing $N$ such that $N\leq (eC^2 h)^{-1} \leq N+1$, which is possible when $h$ is sufficiently small, we can change the right hand side of \r{53} to $\mathcal{O}(e^{-\frac{1}{C^5}})$. 

Now we change variables
\[
(x,y,\eta) \mapsto (x,y,\zeta) = (x,y,\psi_x(y,\eta))
\]
in \r{53} (this is a valid change of variables on a small enough neighborhood of $(0,0,\xi^0)$ by \r{52}). Note that the new phase function, $\Phi(x,y,\zeta)$, satisfies
\be{55}
\begin{split}
\Phi_x(x,x,\zeta) = \zeta, \quad \Phi_y(x,x,\zeta) = -\zeta,\quad \Phi(x,x,\zeta) = 0, \\ 
\mbox{and }\ \mathcal{I} \Phi(x,y,\zeta) > C^4 |x-y|^2 \mbox{ for $x$ near $0$ on the real domain.}
\end{split}
\ee
With this new phase, \r{53} is now in an appropriate form to apply Sjostrand's definition of the analytic wavefront set (see \cite{Sj-A}). Indeed, if there were were only a single amplitude and function in the integral \r{53} we would be finished, as in \cite{FSU}. However, since we have a covector field, $f^s$, and additional function, $\phi$, we must continue a little further. Here we view $[f^s,\phi]$ as a system acted on by an operator with a vector-valued symbol $(A\, b^j, B)$. Note that we have 
\[
\sigma_p(A\, b^j,B)(0,0,\xi^0) = (\tilde a_N(0,\xi^0)\, \tilde b^j(0, \xi^0), \tilde c_N(0,\xi^0)) = (w(0,\theta_0)\, \theta_0, \alpha(0, \theta_0))
\]
using the comments after \r{47} and the fact that $\lambda(0,\theta_0) = 1$ since $\theta_0$ is the velocity vector for a curve in $\Gamma$. Here $\sigma_p$ denotes the principal symbol. Now, following \cite{SU-AJM} we repeat the arguments above for every $(x_0,\theta)$ normal to $(x_0,\xi^0)$ in a neighborhood $W$ of $(x_0,\theta_0)$, which is possible since $U$ is an open subset of $\mathcal{G}$, in order to obtain \r{53} for different symbols and phase functions, with the new symbols satisfying
\be{56}
\sigma_p(A_{\theta}\, b_{\theta}^j,B_{\theta})(0,0,\xi^0) = (w(0, \theta) \, \theta, \alpha(0, \theta)).
\ee
Now let $\{\theta_1, \,...\, , \theta_{n-1}\}$ be a collection of vectors in $W$ which form a basis for the space normal to $\xi^0$. By the elliptic condition there must be some vector $\theta_n = \sum_{j=1}^{n-1} \mu^j \theta_j \in W$ such that 
\be{59}
\G\log  w(x_0,\theta_n) \neq \sum \limits_{j=1}^{n-1} \mu^j \G\log w(x_0, \theta_j).
\ee
This provides us with $n$ equations,
\be{54}
\int_{|x-y|\leq\delta/C^1} e^{\frac{i}{h} \Phi_i(x,\beta)} (A_i (x,\beta;h) \, f^s_j \, b^j_i(x,\beta) + C_i(x,\beta;h) \, \phi) \, \d x = \mathcal{O}\left ( e^{-\frac{1}{hC^5}} \right ), \quad i = 1,\, ... \, ,n
\ee
where for each $i$, $\Phi_i$ is a phase function satisfying \r{55}, and $(A_i\, b_i^j, C_i)$ satisfies \r{56} with $\theta = \theta_i$. We need one more equation to make an elliptic system, and for this we use the fact that $f^s$ is solenoidal on the interior of $M$. Thus $\delta f^s = 0$, and, as in \cite{SU-AJM}, if we take $\chi_0$ to be a smooth cuttoff function with support contained in $M$ and equal to $1$ near $0 \in M$, then $h\, \mathrm{exp}(\frac{i}{h} \Phi_1(x,\beta)) \chi_0 \delta f^s = 0$. Integrating this equality with respect to $x$ yields
\[
0 = \int h\, e^{\frac{i}{h} \Phi_1(x,\beta)} \chi_0(x) \delta f^s(x) \, \d x = - \int e^{\frac{i}{h} \Phi_1(x,\beta)} (i (\Phi_1)_x^j(x,\beta) + h \tilde D^j(x)) f^s_j(x) \, \d x. 
\]
This is equivalent to
\be{57}
\int e^{\frac{i}{h} \Phi_1(x,\beta)} D^j(x,\beta;\sigma) f^s_j(x) \, \d x = 0
\ee
where $\sigma_p(D^j)(0,0,\xi^0) = (\xi^0)^j$. We now claim that the system of $n+1$ equations obtained by combining \r{54} and \r{57} is elliptic at $(0,0,\xi^0)$. By looking at principal symbols we see that this is equivalent to the condition that for a constant covector $f = (f_1, \,...\, ,f_n)$ and constant $\phi$
\[
w(0,\theta_i) (\theta_i)^j f_j  + \alpha(0,\theta_i) \phi = 0 \; i = 1,\,...\, ,n , \quad \mathrm{and}\; (\xi^0)^j f_j = 0
\]
implies $[f,\phi] = 0$. Indeed this is true by the ellipticity condition, \r{59}, and the fact that the $\theta_i$ contain a basis for the space perpendicular to $\xi^0$. Let us now rewrite the system in the more compact form
\be{101}
\int_{|x-y| \leq C} \mathrm{diag}(e^{\frac{i}{h} \Phi_i(x,\beta)}) \mathbf{A}(x,\beta) [f^s,\phi](x) \, \d x = \mathcal{O}(e^{-\frac{1}{hC}})
\ee
where $\mathbf{A}(x,\beta)$ is a matrix valued symbol acting on $[f^s,\phi]$ and $C$ is a new constant. The ellipticity of the system at $(0,0,\xi^0)$ is simply equivalent to the invertibility of the principle symbol of $\mathbf{A}$ at this point.

We would like to replace $\mathbf{A}$ in \r{101} by the identity matrix so that, with this replacement, \r{101} precisely shows, in the sense of \cite[Def.~6.1]{Sj-A}, that $(x_0, \xi^0) \notin WF_A([f^s,\phi])$. To do this we generalize the proof of proposition 6.2 from \cite{Sj-A} to the case of operators with matrix valued symbols. Indeed, following \cite{Sj-A} and \cite{SU-AJM} we first define a system of \PDO's in the complex domain by writing
\be{102}
\mathrm{Op}(\mathbf{A}) [f^s,\phi] (y) = \iint \mathrm{diag}(e^{\frac{i}{h} (\Phi_i(y,\beta) - \overline{\Phi_i(x,\beta)}}) \mathbf{A}(x,\beta) [f^s,\phi](x) \, \d x \, \d \beta.
\ee
These operators have differing phase functions, but by making appropriate changes of variables we can change them all to the same phase $\Phi$ and this does not change the principal symbol. Therefore, we may construct a parametrix for the system and following \cite{Sj-A} use this parametrix to express $\mathbf{Id} \,e^{\frac{i}{h} \Phi}$ as a superposition of functions $\mathbf{A} \,e^{\frac{i}{h} \Phi}$ modulo an exponentially decreasing function. Following now the same argument as is given for proposition 6.2 in \cite{Sj-A}, but with matrix valued symbols, we have (possibly with a new constant $C$)
\be{103}
\int_{|x-y| \leq C} e^{\frac{i}{h} \Phi(x,\beta)} \mathbf{Id} \,[f^s,\phi](x) \, \d x = \mathcal{O}(e^{-\frac{1}{hC}}),
\ee
for $\beta = (y,\eta)$ in a neighborhood of $(0,\xi^0)$. This proves that $(0,\xi^0)$ is not in $\mathrm{WF_A}([f^s,\phi])$ as originally claimed.
\end{proof}

\begin{proof}[Proof of Theorem~\ref{thm1}]

Let $f \in L^2(M)$ be such that $I_{G,\lambda,w}f = 0$, and pick an intermediate manifold $M_{1/2}$ with $M \Subset M_{1/2} \Subset M_1$. Now extend $f$ as zero outside of $M$ and consider the decomposition \r{31} on $M_{1/2}$
\be{58}
f = f^s_{M_{1/2}} + \d \phi_{M_{1/2}}.
\ee
We extend $f^s_{M_{1/2}}$ and $\phi_{M_{1/2}}$ as zero outside of $M_{1/2}$ and then apply Proposition~\ref{pr3} to the pair $[f^s_{M_{1/2}},\phi_{M_{1/2}}] \in \mathbf{L}^2(M_1)$ with $M$ replaced by $M_{1/2}$. The proposition implies that the pair $[f^s_{M_{1/2}},\phi_{M_{1/2}}]$ is analytic in the interior of $M_{1/2}$, and so by \r{58} $f$ is also analytic in the interior of $M_{1/2}$. Since $f$ is identically zero on $M_{1/2}\setminus M$ this implies that $f$ must be equal to zero on all of $M$. Therefore the kernel of $I_{G,\lambda,w}$ is just $\{0\}$, and so $I_{G,\lambda,w}$ is injective.

\end{proof}

\section{Application to the attenuated ray transform of vector fields}  \label{sec4}
The Riemannian version of the attenuated ray transform of vector fields, see \r{02}, \r{03}, is given by the weighted transform \r{04} with a weight
\be{P1}
w_\sigma (x,\theta) = e^{-\int_{\tau_-(x,\theta)}^0 \sigma(\gamma_{x,\theta}(s), \dot \gamma_{x,\theta}(s)) \, \d s}.
\ee
Compare this with \r{03}, where $\sigma(x,\theta)$ is a given, usually non-negative, function on $SM$ modeling the attenuation properties of the media, and the generator $\mathbf{G}$ is the generator of the geodesic flow. For the Riemannian case of the attenuated ray transform and its relation to the transport equation, we refer to \cite{Sh-transport}. In the analysis that follows, actually $\mathbf{G}$ does not need to be the  geodesic  generator. The weight \r{P1} satisfies the relation
\[
\mathbf{G}\log w = -\sigma.
\]
The results in section~\ref{sec2} then imply the following.
\begin{corollary}
The elliptic condition is satisfied in any of the following cases:

(a) $w=w_\sigma$ is with $\sigma\not=0$ on $\mathcal{G}$.

(b) $U$  satisfies \r{U} and $w=\chi w_\sigma$, where $\mathbf{G}\chi=0$, $\chi|_U\not=0$, and 
 either $\sigma=\sigma(x)\not=0$ is a function of $x$ only, or $\sigma=\sigma_{ij}(x)\theta^i\theta^j\not=0$.
\end{corollary}

In case (a), we integrate over all geodesics and we require the absorption $\sigma$ to be non-zero everywhere, although possibly anisotropic (dependent on $\theta$ as well as $x$). In case (b), we restrict the ray transform to  $U$ using the smooth cut off $\chi$ but we then require $\sigma$ to be isotropic and non-zero; or a non-degenerate quadratic form of $\theta$ for any $x$. One can allow more general $\sigma$'s in (b) but we only mention cases that might be useful in applications. 

\begin{proof}
If $w$ is as in (a), then we can choose $U = \mathcal{G}$. Then $\mathbf{G} \log w=-\sigma$ cannot be of the form \r{econd} on $\xi^\perp$ for a fixed $x$ because it never vanishes. 

Assume (b). Then either $\mathbf{G}\log w= - \sigma(x)$ or
$\mathbf{G}\log w=- \sigma_{ij}(x)\theta^i\theta^j$ 
on the set $\chi\not=0$; and in particular in $U$. In either case, it is not of the form \r{econd} on $\xi^\perp \cap U$. 
\end{proof}

\textbf{Acknowledgements.} The authors thank J.~Sj\"ostrand for his help with certain aspects of the analytic microlocal theory presented in \cite{Sj-A}.

\newpage
\bibliographystyle{abbrv}
\bibliography{Doppler.bbl}

\end{document}